\documentclass[]{amsart}
\usepackage{amscd}
\usepackage{amssymb}
\usepackage{amsfonts}
\usepackage{remark}
\usepackage{enumerate}

\input BoxedEPS.tex
\SetRokickiEPSFSpecial
\HideDisplacementBoxes

\newtheorem{thm}{Theorem}[section]
\newtheorem{theorem}[thm]{Theorem}
\newtheorem{lemma}[thm]{Lemma}
\newtheorem{corollary}[thm]{Corollary}
\newtheorem{proposition}[thm]{Proposition}
\newremark{definition}[thm]{Definition}
\newremark{remark}[thm]{Remark}
\newremark{emp}[thm]{\kern -4pt}

\newcommand{\W}{{\mathcal W}}
\newcommand{\X}{{\mathcal X}}
\newcommand{\Y}{{\mathcal Y}}
\newcommand{\Z}{{\mathcal Z}}

\newcommand{\cO}{\mathcal O}
\newcommand{\Spec}{{\rm Spec}\kern 2pt}
\newcommand{\Proj}{{\rm Proj}\kern 2pt}
\newcommand{\Gal}{{\rm Gal}\kern 1pt}
\newcommand{\Aut}{{\rm Aut}\kern 1pt}
 
\newcommand{\wt}{\widetilde}

\begin{document}

\title{\bf Stable reduction of finite covers of curves}

\author{Qing Liu}


\begin{abstract} Let $K$ be the function field of a connected 
regular scheme $S$ of dimension $1$, and let $f : X\to Y$ be
a finite cover of projective smooth and geometrically connected
curves over $K$ with $g(X)\ge 2$. Suppose that $f$ can be 
extended to a finite cover $\X \to \Y$ of semi-stable models 
over $S$ (it is known that this is always possible up to finite 
separable extension of $K$).  Then there exists a unique minimal 
such cover. This gives a canonical way to extend $X\to Y$ to a 
finite cover of semi-stable models over $S$. 
\end{abstract} 

\address{CNRS, Laboratoire A2X, Institut de Math\'ematiques de Bordeaux\\ 
Universit\'{e} de Bordeaux I\\351, cours de la Lib\'{e}ration, 33405 
Talence, FRANCE}

\maketitle

Let $S$ be a Dedekind scheme (i.e a connected Noetherian regular scheme 
of dimension $1$), with field of functions $K:=K(S)$. Let $f : X\to Y$ 
be a finite morphism of smooth geometrically connected 
projective curves over $K$. We can ask how to extend,
in some canonical way, the morphism $f$ to a morphism of models 
of $X$ and $Y$ over $S$. It 
is proved in \cite{LL},4.4 that if $X$ and $Y$ have stable models 
$\X^{\rm st}, \Y^{\rm st}$ over $S$, then $f$ extends uniquely to a 
morphism $\X^{\rm st}\to \Y^{\rm st}$. However we will in general
lose the finiteness of $f$. On the other hand, after a finite 
separable extension of $K$, $f$ extends to a finite morphism of 
semi-stable models $\X\to\Y$ over $S$
(\cite{Col}; \cite{LL}, Remark 4.6; and Corollary \ref{sst.red} below).
Following Coleman \cite{Col}, such a pair $\X\to\Y$ is called a 
{\em semi-stable model} of $f$, and it is called {\em stable} if 
moreover it is minimal among the semi-stable models of $f$ 
(cf. \ref{model-f}). The stable model of $f$ (if it exists) 
is unique up to isomorphism.

\begin{theorem} \label{main-cor} 
{\rm (Corollary \ref{stable})} 
Suppose that either $g(X)\ge 2$, or $g(X)=1$ and $X$ has 
potentially good reduction. Then there exists a finite separable 
extension $K'$ of $K$ such that $X_{K'}\to Y_{K'}$ admits a stable model 
$\X'\to\Y'$ over $S'$, where $S'$ is the integral closure of $S$
in $K'$. Moreover, for any Dedekind scheme $T$ dominating $S'$, 
$\X'_{T}\to \Y'_{T}$ is the stable model of $X_{K(T)}\to Y_{K(T)}$
over $T$. 
\end{theorem} 

This gives a canonical way to extend a finite cover of projective 
smooth curves over $K$ to a finite cover of semi-stable models
over (some finite cover of) $S$. The last part of the
theorem says that the stable model of $f$ commutes with flat base
change. The stable model of $f$ should be seen as an analogue
of the stable model of a curve. 
In a forthcoming work, we will apply this theorem to study 
a compactification of Hurwitz moduli spaces of finite covers of 
curves. We also prove the theorem for smooth marked curves $X, Y$.
Note that in general $\X', \Y'$ are not the respective stable 
models of the curves $X_{K'}, Y_{K'}$. 
The proof of \ref{main-cor} is based on a more 
general result: 

\begin{theorem} \label{main} {\rm (Theorem \ref{stable-e})} 
Let $f : X\to Y$ be a finite morphism of 
smooth geometrically connected projective curves over $K$. Let $\X, \Y$ 
be respective models of $X$ and $Y$ over $S$. Then there exists
a finite separable extension $K'$ of $K$, and models $\X', \Y'$ of $X_{K'}$ 
and $Y_{K'}$ over $S'$ (integral closure of $S$ in $K'$) 
such that 
\begin{enumerate}[\hskip 8pt \rm (1)] 
\item $\X', \Y'$ dominate respectively $\X_{S'}, \Y_{S'}$ and are 
semi-stable over $S'$; 
\item $X_{K'}\to Y_{K'}$ extends to a finite morphism $\X' \to \Y'$;
\item the pair $(\X', \Y')$ is minimal : for any pair $(\X'', \Y'')$
satisfying properties {\rm (1)} and {\rm (2)}, $\X''$ and $\Y''$
dominate respectively $\X'$ and $\Y'$. The cover $\X'\to \Y'$ is called
the {\em stable hull of} $\X_{S'}\dasharrow\Y_{S'}$. 
\item the formation of $\X', \Y'$ commutes with flat base change
$T \to S'$: if $T$ is a Dedekind scheme dominating $S'$, then 
$\X'_T\to \Y'_T$ is the stable hull of $\X_T\dasharrow\Y_T$. 
\end{enumerate}
\end{theorem} 

The paper is organized as follows. In Section 1, $S=\Spec\cO_{K}$ 
is local. We discuss some sufficient conditions for a model 
of $X_{\hat{K}}$ over $\hat{\cO}_K$ (completion of $\cO_K$) 
to be defined over $\cO_K$. 
E.g., this is true for semi-stable or regular models. 
This is used to reduce the proof of \ref{finite-h} to the case 
of a complete local base $S$. In Section 2, we prove 
the special case of \ref{main} when $X=Y, \X=\Y$ and
$f={\rm Id}$. Parts (1) and (2) of \ref{main} are proved
in Section 3. The existence of a stable model (Part (3))
is proved in Section 4. 
\medskip

{\noindent\bf Convention} 
Through this paper, $S$ is a Noetherian regular connected scheme
of dimension $1$ (Dedekind scheme), and, unless otherwise 
specified, $X, Y$ are projective, smooth and geometrically 
connected curves over $K:=K(S)$. When we state that a property 
(P) holds for some model $\X$ of $X$ {\em after finite separable 
extension of $K$}, this means that there exists a finite separable
extension of $K'/K$, such that (P) is satisfied for 
$\X_{S'}$, where $S'$ is the normalization of $S$ in $K'$. 

Note that the hypothesis $X, Y$ geometrically connected 
(instead of connected) is not serious. Actually $X$ is
geometrically connected over the finite separable 
extension $H^0(X,\cO_X)$ of $K$. 
\medskip 


\noindent {\sc Acknowledgements}
\smallskip 

I would like to thank D. Lorenzini for pointing out Coleman's
paper \cite{Col} and for interesting comments, Ch. Deninger 
and A. Werner for motivational 
discussions and Sylvain Maugeais for pointing out some mistakes
in an earlier version. I also thank the referee to correct a
mistake in the proof of Proposition \ref{aut-stable}. 


\begin{section}{Descent from the completion} 

Let $\cO_K$ be a discrete valuation ring, and let $X$ be 
a geometrically connected smooth projective curve $K$. 
We give some sufficient conditions for a model of $X_{\hat{K}}$
over $\hat{\cO}_K$ to be defined over $\cO_K$.  

\begin{lemma} \label{descent}
Let $\cO_K$ be a discrete valuation ring, let $A$ be a flat 
local $\cO_K$-algebra, localization of a finitely generated 
$\cO_K$-algebra. Suppose that 
$\Spec (A\otimes {\hat{K}})$ is integral, one-dimensional, and smooth
over $\hat{K}$. 
Let $\W\to \Spec \hat{A}$ be a projective birational morphism. 
Then under any of the following conditions,  
$\W$ and the morphism ${\W}\to\Spec{\hat{A}}$ are defined over $A$:
\begin{enumerate}[{\hskip 8pt\rm (a)}]
\item ${\rm Pic}(\hat{A}\otimes {\hat{K}})$ is a torsion group 
(i.e., every element has finite order);
\item $\hat{A}$ is $\mathbb Q$-factorial. 
\end{enumerate}
\end{lemma} 

{\noindent\it Proof:} The morphism $\W\to \Spec(\hat{A})$ is the  
blowing-up along a closed subscheme $V(\mathfrak {I})$ of $\Spec(\hat{A})$. 
Let $t$ be a uniformizing element of $\cO_K$. 
Let us show that $\mathfrak I$ can be chosen in such a way that 
$t^n \in \mathfrak I$ 
for some $n\ge 1$. Let us suppose that condition (a) is 
satisfied. Then there exists some $m\ge 1$ such that the restriction 
of $\mathfrak I^m$ to the generic fiber is principal, generated by a 
$f\in \mathfrak I^m\otimes \hat{K}$. Under condition (b), 
let $D$ be the scheme theoretical closure of 
$V(\mathfrak I)\cap \Spec(\hat{A}\otimes \hat{K})$
in $\Spec(\hat{A})$, considered as a Weil divisor. Then $mD$ is principal for
some $m>0$, generated by a $f\in \mathfrak I^m$. 

In both cases, there exist $a, b\in \mathbb Z$ such that 
$t^a\mathfrak I^m\subseteq f\hat{A} \subseteq t^b \mathfrak I^m$. 
Replacing $\mathfrak I$ by the ideal $t^{a}f^{-1}\mathfrak I^m$ 
(which does not change the blowing-up along $\mathfrak I$), 
we can suppose that the restriction of $\mathfrak I$ 
to the generic fiber is trivial, and $t^n\in \mathfrak I$ for some $n\ge 0$
($n=a-b$). Since $\cO_K$ is dense in $\hat{\cO}_K$, 
$\mathfrak I$ is then generated by an ideal $I$ of $A$. Let 
$\X\to\Spec A$ be the blowing-up along $V(I)$, then 
$\W \to \Spec(\hat{A})$ is obtained from $\X\to\Spec A$ 
by the base change $\hat{\cO}_K/\cO_K$.
\qed 

\begin{remark} \label{reg-d} 
Let $\W$ be a projective regular model of $X_{\hat{K}}$ over 
${\hat{\cO}_K}$ dominating $\X_{\hat{\cO}_K}$ for some
model $\X$ of $X$ over $\cO_K$. 
Then $\W$ and $\W\to \X_{\hat{\cO}_K}$
are defined over $\cO_K$. Actually, if $\wt{\X}\to\X$ is
the minimal desingularization of $\X$, then 
$\wt{\X}_{\hat{\cO}_K}\to\X_{\hat{\cO}_K}$ is the minimal
desingularization of $\X_{\hat{\cO}_K}$ 
(same proof as in \cite{Liub}, 9.3.28). 
Hence $\W\dasharrow \wt{\X}_{\hat{\cO}_K}$ is a birational 
morphism. It is a sequence of blowing-ups of closed points
(\cite{Liub}, 9.2.2). 
Hence $\W$ is defined over $\cO_K$, and so is 
the morphism $\W\to\X_{\hat{\cO}_K}$ by faithfully flat descent
of rational maps.
\end{remark} 

\begin{proposition} \label{completion} Let $S=\Spec\cO_{K}$ be local. 
Let $\X$ be a model of $X$ over $S$ dominating some semi-stable or
regular model. Let $\varphi : \W\to \X_{\hat{\cO}_K}$ be a projective 
birational morphism over $\hat{\cO}_K$. Then 
$\W$ and $\varphi$ are defined over $S$. 
\end{proposition} 

{\noindent\it Proof:} It is enough to show that $\W$ is defined over 
$S$. We can suppose that $\X$ itself is semi-stable or regular. 
The morphism $\varphi$ is an 
isomorphism outside of a finite set $F$ of closed points of 
$(\X_{\hat{\cO}_K})_s=\X_s$. Let $x\in F$, 
and let $A_x=\cO_{\X,x}$. If $\X$ is regular, then $\hat{A}_x$ is
regular and thus factorial. If $\X$ is semi-stable, 
there exists a finite (\'etale) extension $\cO_L/\hat{\cO}_K$ such that
$\hat{A}_x\otimes_{\hat{\cO}_K}{\cO_L}$ is a finite direct sum of 
rings $\cO_L[[u]]$ or $\cO_L[[u,v]]/(uv-a)$, $a\in {\cO}_L$. 
It is well-known that ${\rm Pic}$ of these rings tensored by 
$\hat{L}$ are trivial (see for instance \cite{He}, Cor. 2.2, for 
${\cO}_L[[u,v]]\otimes \hat{L}/(uv-a)$).
Hence ${\rm Pic}(\hat{A}_x)$ is torsion 
(\cite{Liub}, 7.2.18, 7.2.19). By Lemma \ref{descent}, 
$\W\times_{\X_{\hat{\cO}_K}}
\Spec\hat{A}_x\to\Spec\hat{A}_x$ is defined over $A_x$.
By glueing the morphisms above $\Spec A_x$, when $x$ varies in $F$, 
and the isomorphism above $\X\setminus F$, we find a morphism over $\X$
which is equal to $\varphi$ when base changed to $\hat{\cO}_K$.
\qed 

\begin{remark} In general, not every birational projective morphism
$\W\to \X_{\hat{\cO}_K}$ is defined over $\cO_K$. In other words,
even if $\W\to \X_{\hat{\cO}_K}$ is an isomorphism outside of the
special fiber, it is not necessarily the blowing-up along a closed
subscheme with support in the special fiber. To construct such a
counterexample, we will imitate the example of a non contractible
component given in \cite{BLR}, Lemma 6.7.6. 
Let us consider the smooth elliptic curve $\mathcal E/\cO_K$ 
and the point $a_k\in {\mathcal E}_k(k)$ as given in \cite{BLR}, p. 171. 
The point $a_k$ satisfies the property that no multiple (in the sense of the
group law on $\mathcal E$) $na_k$, $n>0$, can be lifted to a section in 
$\mathcal E(\cO_K)$. 

\if false
\begin{figure}
\centering
\includegraphics{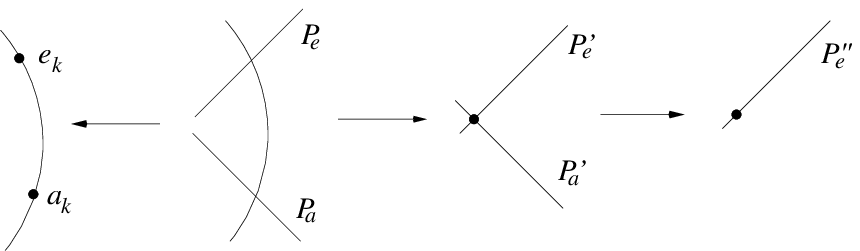}
$$ \mathcal E_k \qquad \longleftarrow\qquad
\X'_k \qquad\longrightarrow\qquad \W_k 
\qquad\longrightarrow\qquad \X_k$$
\end{figure} 
\fi 

\newpage 
\medskip
\centerline{\BoxedEPSF{contract.eps}}
\medskip

$$ \mathcal E_k \qquad \longleftarrow\qquad
\X'_k \qquad\longrightarrow\qquad \W_k 
\qquad\longrightarrow\qquad \X_k$$ 

Let $\X'\to \mathcal E$ be the blowing-up of $\mathcal E$ along $\{a_k, e_k\}$,
where $e$ is the unit section of $\mathcal E$. Let $P_a$ and $P_e$ be the
respective inverse images of $a_k, e_k$ in $\X'$.  They are projective
lines over $k$. Let $\wt{\mathcal E}_k$ be the strict transform
of ${\mathcal E}_k$ in $\X'$. The unit section of $E_K$ gives 
rise to a section of $\X'$ which meets $\X'_k$ at an interior point 
of $P_e$. Hence there exists a contraction map $\X'\to \X$ of 
$P_a\cup \wt{\mathcal E}_k$ (\cite{BLR}, Corollary 6.7.3). Now over 
$\hat{\cO}_K$, we can contract $\wt{\mathcal E}_k$ in $\X'_{\hat{\cO}_K}$ 
({op.\kern 2pt cit.,} 6.7.4), which gives 
a model $\W/{\hat{\cO}_K}$ and we have birational projective morphisms
$\X'_{\hat{\cO}_K}\to \W\to \X_{\hat{\cO}_K}$. We want to show that 
$\W$ is not defined over $\cO_K$, or, equivalently, that $\wt{\mathcal E}_k$
can not be contracted over $\cO_K$. Suppose that $\W$ exists over $\cO_K$,
then there exists a relative effective Cartier divisor $D$ on $\X'$ such that
${\rm Supp} D$ meets $P_a$ and $P_e$ but not $\wt{\mathcal E}_k$. Taking the
Zariski closure of $D_K$ in $\mathcal E$ gives a relative effective Cartier 
divisor $D'$ on $\mathcal E$ such that $D'_k=na_k+me_k$ for some integers 
$n, m>0$. 
But $D'-(n+m)e\in {\rm Pic}^0_{\mathcal E/\cO_K}(\cO_K)\simeq \mathcal E(\cO_K)$ 
then defines a section $b\in\mathcal E(\cO_K)$ such that $D'-(n+m)e\sim 
b - e$. In the special fiber we then have $n(a_k-e_k)\sim b_k-e_k$,
hence $na_k=b_k$ in the group $\mathcal E_k(k)$, 
contradiction with the assumption on $a_k$.
\end{remark} 

\begin{lemma} \label{rm-d} 
Suppose that $X$ has semi-stable reduction over $S$. Then 
every relatively minimal semi-stable model $\W$ of $X_{\hat{K}}$ 
over $\hat{\cO}_K$ is defined over $\cO_K$.
\end{lemma} 

{\noindent\it Proof:} 
If $g(X)\ge 2$, then the unique relatively minimal semi-stable
model of $X$ is the stable model of $X$ over $S$. 
Since the stable model is unique and commutes with flat
base change (\cite{Liub}, 10.3.36), $\W$ 
over $\hat{\cO}_K$ is defined over $\cO_K$. 
The same is true if $g(X)=1$ and $X$ has good reduction.   

Suppose that $g(X)=1$ and $X$ has multiplicative reduction. 
Let $\X$ be a relatively minimal semi-stable model. Let
$\wt{\X}$ be its minimal desingularization. Let us show that
$\wt{\X}$ is the minimal regular model of $X$ over $S$.
Let $\Gamma$ be an irreducible component of $\X_s$ and let $\wt{\Gamma}$
be its strict transform in $\wt{\X}$. By Lemma \ref{contract-2}(a),
$\deg\omega_{\X/S}|_{\Gamma} > 0$ if $\X_s$ is reducible. 
If $\X_s$ is irreducible, then $\deg\omega_{\X/S}|_{\Gamma}=2g(X)-2=0$. 
So in any case 
$\deg\omega_{\wt{\X}/S}|_{\wt{\Gamma}}=\deg\omega_{\X/S}|_{\Gamma} \ge 0$. 
Hence $\wt{\X}$ has no exceptional divisor. Since 
any strict subset of the set of irreducible components of $\wt{\X}_s$
can be contracted (\cite{Liub}, 9.4.19) into a semi-stable model
(\ref{contract-2}(a)), $\X_s$ is irreducible. 
So the relatively minimal semi-stable models of $X$ correspond
bijectively to the irreducible components of the minimal regular
model of $X$. The same is true over $\hat{\cO}_K$. 
Since the minimal regular model commutes with
the base change $\hat{\cO}_K/\cO_K$ (\cite{Liub}, 9.3.28), we 
see that the relatively minimal semi-stable models of 
$X_{\hat{K}}$ are exactly those of $X$ base changed to $\hat{\cO}_K$. 

Suppose that $g(X)=0$. Let $\X$ be a relatively minimal semi-stable
model of $X$ over $S$ and let $\Gamma_1,...,\Gamma_n$ be the irreducible
components of ${\X}_s$. Then 
$$\sum_i \deg\omega_{\X/S}|_{\Gamma_i}=2g(X)-2 <0.$$ 
So $\omega_{\X/S}$ has 
negative degree on at least one component $\Gamma_i$. By 
Lemma \ref{contract-2}(a), $n=1$, thus $\X_s$ is irreducible,
semi-stable and has arithmetic genus $0$. So $\X_s$ is smooth. 
The same reasoning shows that $\W$ is smooth. 
Let $\mathcal L$ be the dual of the dualizing sheaf on 
$\W$. Let $s_0, s_1\in H^0(\W, \mathcal L)$ be a basis
with $s_i\in H^0(X, \omega_{X/K}^{\vee})$. 
Then $s_0, s_1$ define a closed immersion $\W\to \mathbb P^2$ 
over $\hat{\cO}_K$. Its image is a conic defined by a polynomial 
with coefficients in $K$. Hence $\W$ is defined over $\cO_K$.
\qed 
\smallskip 

The next corollary is weaker than \ref{completion}, but 
sufficient for the purpose of Theorem \ref{stable-e}. 
It is an immediate consequence of \ref{completion} and the above lemma.
However we give a direct proof without using \ref{completion}. 

\begin{corollary} \label{st-d} Let $S=\Spec\cO_K$ be local.
Let $X$ be a smooth geometrically connected projective curve over $K$. Then
every semi-stable model $\W$ of $X_{\hat{K}}$ over $\hat{\cO}_K$ 
is defined over $\cO_K$. 
\end{corollary} 

{\noindent\it Proof:} The model $\W$ dominates a relatively
minimal semi-stable model ${\Z}$ of $X_{\hat{K}}$ which is
then defined over $\cO_K$ by \ref{rm-d}. Let $\wt{\W}\to\W$ 
be the minimal desingularization of $\W$.
Then $\wt{\W}\to {\Z}$ is defined over $\cO_K$ 
(\ref{reg-d}). The irreducible components of $\wt{\W}_s$ 
in the exceptional locus of $\wt{\W}\to\W$ are $(-2)$-curves 
and can be contracted over $\cO_K$ (\ref{contract-2}(a)). Hence
$\W$ is defined over $\cO_K$. 
\qed 
\smallskip 

\begin{remark}
If $X=\mathbb P^1_K$, then every normal model $\W$ 
of $X_{\hat{K}}$ is defined over $\cO_K$. Indeed, if $\wt{\W}$ 
is the minimal desingularization of $\W$, then
$\wt{\W}$ dominates a relatively minimal regular model. The latter 
is smooth and defined over $\cO_K$. Hence $\wt{\W}$ is defined over
$\cO_K$. Now every strict subset of the set of irreducible
components of $\wt{\W}_s$ is contractible over $\cO_K$ 
(\cite{Liub}, Exerc. 9.4.5). 
So $\W$ is defined over $\cO_K$. 
\end{remark}

\end{section}

\begin{section}{Stable hull of a model}

\begin{definition} Let $S$ be a connected Noetherian regular scheme
of dimension $1$ (i.e., a {\em Dedekind scheme}). 
Let $X$ be an integral projective variety over 
$K$. A {\em model $\X$ of $X$ over $S$} is an integral projective 
scheme over $S$ whose generic fiber is isomorphic to $X$. 
Recall that $\X$ is said to be {\em semi-stable} if its geometric
fibers are reduced with only ordinary double points as singularities. 
A morphism
of models is defined in an obvious way. 
\end{definition}  

\begin{definition} Let $X$ be a connected projective smooth 
curve over $K$, and let $\X$ be a model of $X$ over $S$. A {\em stable hull} 
of $\X$ is a semi-stable model $\W$ of $X$ dominating $\X$, and 
minimal for these properties (i.e., every semi-stable model dominating
$\X$ dominates $\W$).  
\end{definition} 

The aim of this section is to prove the next result.

\begin{thm} \label{ss-h} Let $X$ be a geometrically 
connected projective smooth curve 
over $K$ and let $\X$ be a model of $X$ over $S$. 
\begin{enumerate}[\hskip 8pt \rm (a)] 
\item The stable hull of $\X$ is unique (up to isomorphism) 
when it exists. In general, there 
exists a finite separable extension $K'/K$ such that 
$\X_{S'}$ (where $S'$ is the integral closure of $S$ in $K'$) 
has a stable hull over $S'$. 
\item The stable hull commutes with flat base change: suppose that 
$\X$ admits a stable hull $\W$ over $S$ and let $S'\to S$ be a flat 
morphism of Dedekind schemes, then $\W_{S'}$ is the stable 
hull of $\X_{S'}$ over $S'$.
\end{enumerate} 
\end{thm} 

The proof of the theorem is postponed to \ref{pf-ssh}. 

\begin{lemma} \label{equiv-lm} 
Let $G$ be a finite group acting on $X$. Let $\X$ be a model of 
$X$ over $S$. Then there exists a model 
$\Z$ of $X$ dominating $\X$, endowed with an action of $G$, and
minimal for these properties. 
\end{lemma} 

{\noindent\it Proof:} (See \cite{dJ0}, 7.6) 
Let $\sigma\in G$, then there exists a 
model $\X^{\sigma}$ such that $\sigma : X\to X$ extends to 
an isomorphism $\sigma : \X\to \X^{\sigma}$. If $\tau \in G$, 
by composing $(\sigma^{-1}\tau) : \X^{\tau^{-1}\sigma}\to \X$ with
$\sigma : \X\to \X^{\sigma}$, 
we obtain an isomorphism $\X^{\tau^{-1}\sigma}\to \X^{\sigma}$ 
denoted by $\tau$. Let $\mathcal P$ be the fiber 
product $\prod_{S, \sigma\in G} \X^{\sigma}$ over $S$. Then 
we can make $G$ act on $\mathcal P$ by 
$$\tau : (x_{\sigma})_{\sigma} \mapsto (\tau(x_{\tau^{-1}\sigma}))_{\sigma}.$$
Moreover, the diagonal morphism $\Delta : X\to \mathcal P_K, x\mapsto 
(x,...,x)$ is $G$-equivariant. Let $\Z$ be the Zariski closure of 
$\Delta(X)$ in $\mathcal P$, endowed with the reduced structure.
Then $G$ leaves stable $\Z$. Note that $\Z$ dominates $\X$ because
the projection morphism $\mathcal P\to \X$ induces a morphism 
$\Z\to \X$ which is an isomorphism on the generic fibers. 

Let us prove that $\Z$ is minimal. Let $\W$ be a model endowed with
an action of $G$ and a birational morphism $\W\to \X$. Then we have
a morphism $\W\to\X^{\sigma}$ for all $\sigma$ and hence a morphism
$h : \W\to\mathcal P$. 
Since $h(\W)$ is irreducible and contains $Z$, $h$ induces a morphism 
of models $\W\to\Z$. 
\qed  
\medskip 

Let  us give a corollary of \ref{ss-h}. 

\begin{corollary} \label{equiv} Let $G$ be a finite group acting on $X$. Let
$\X$ be a model of $X$ over $S$. Then after finite separable
extension of $K$, $\X$ is dominated by a semi-stable (resp. semi-stable
and regular) model $\W$ such that the action of $G$ extends to $\W$. 
Moreover, there exists a minimal such a model $\W$. 
\end{corollary} 

{\noindent\it Proof:} Let $\Z$ be the model defined in \ref{equiv-lm}.
Let $\W$ be the stable hull of $\Z_{S'}$ over some $S'/S$
(Theorem \ref{ss-h}(a)). By the uniqueness property, the action of 
$G$ on $\Z_{S'}$ extends to $\W$. It is clear that $\W$ is minimal
with respect to the required properties. To have a minimal semi-stable
regular model, it is enough to take the minimal desingularization of 
$\W$. 
\qed 

\begin{lemma} \label{reduced-f} Let $S$ be local with separably
closed residue field. Let $\X$ be a model of $X$ with minimal 
desingularization $\wt{\X}\to\X$. Suppose that $\wt{\X}$ dominates 
a regular model $\Z$ 
and that $\X_s$, $\Z_s$ are geometrically reduced. Then $\wt{\X}_s$
is geometrically reduced.
\end{lemma} 

{\noindent\it Proof:} Note that if $\W$ is a normal model of $X$,
then $\W_s$ verifies the property (S$_1$), thus $\W_s$ is geometrically 
reduced if and only if every irreducible 
component of $\W_s$ has geometric multiplicity (\cite{BLR}, 9.1.3)
equal to $1$ in $\W_s$. The latter condition depends only on the
generic points of $\W_s$. We can decompose $\wt{\X}\to \Z$ into 
a sequence of blowing-ups
$$ \wt{\X}=:\Z_0\to \Z_{1} \to ... \to \Z_n:=\Z$$
such that $\Z_{i} \to \Z_{i+1}$ consists in blowing-down
an exceptional divisor $\Theta_i$ contained in $\Z_i$. 
We will show by induction that $\Theta_i$ has multiplicity $1$ in 
$(\Z_i)_s$ and $\Theta_i\simeq\mathbb P^1_{k(s)}$.
Since $(\Z_n)_s$ is geometrically reduced, this will imply that
$\wt{\X}_s$ is geometrically reduced. 

By the minimality of $\wt{\X}\to\X$, $\Theta_0$ is not 
mapped to a closed point in $\X$. Thus $\wt{\X}\to\X$ is
an isomorphism in a neighborhood of the generic point
of $\Theta_0$. In particular, $\Theta_0$ has geometric multiplicity
$1$. Since 
$\Theta_0$ is a projective line by Castelnuovo criterion, 
it is isomorphic to $\mathbb P^1_{k(s)}$. 
Suppose that the same holds for $\Theta_{j}$, $0\le j\le i-1$. 
Let ${\Theta}_{i,j}$ be the strict transform of 
$\Theta_i$ in $\Z_j$, $0\le j\le i-1$. If $\Theta_{i,j}$ meets
$\Theta_{j}$ for some $j$, then the computation of $\Theta_j^{2}=-1$
shows that ${\Theta}_{i,j}$ has multiplicity $1$ and cuts 
$\Theta_{j}$ at a rational point. Thus $\Theta_i$ has 
multiplicity $1$ and is isomorphic to $\mathbb P^1_{k(s)}$. 
Otherwise, $\Z_0\to \Z_i$ is an isomorphism in a neighborhood of 
$\Theta_{i,0}$. In particular, $\Theta_{i,0}$ is an exceptional
divisor. Then we can conclude exactly as for $\Theta_0$.
\qed 

\begin{proposition} \label{st-hull.ex}
Let $\X$ be a normal model of $X$ over $S$. Then the following 
properties are equivalent:
\begin{enumerate} [{\rm \hskip 8pt (i)}]
\item $\X$ is dominated by a semi-stable model over $S$;
\item $X$ admits a semi-stable model over $S$ and $\X_{{s}}$ 
is geometrically reduced for all $s\in S$; 
\item  The minimal desingularization $\wt{\X}$ of $\X$ is 
semi-stable over $S$. 
\end{enumerate}
\end{proposition} 

{\noindent\it Proof:}
(i) $\Longrightarrow$ (ii) If $\X$ is 
dominated by a semi-stable model $\X'$, then any irreducible 
component $\Gamma$ of $\X_s$ is birational to an irreducible 
component of $\X'_s$. The latter being geometrically reduced, 
$\X_{{s}}$ is geometrically reduced.

(ii) $\Longrightarrow$ (iii) We know that $\wt{\X}$ dominates a relatively 
minimal regular model $\Z$. The semi-stable reduction hypothesis 
implies that $\Z$ is semi-stable (see \cite{Liub}, 10.3.34(a) if
$g(X)\ge 1$; if $g(X)\le 0$, then $\Z$ is smooth). 
Since the minimal desingularization commutes with
\'etale base change (see the proof of \cite{Liub}, 
Prop. 9.3.28), we can suppose that $S$ is local 
with separably closed residue field. By Lemma \ref{reduced-f}, 
$\wt{\X}_s$ is geometrically reduced. 
The map $\wt{\X} \to \Z$ consists in blowing-up successively 
closed points. The fact that $\wt{\X}_s$ is geometrically
reduced implies that we only blow-up rational points in the
smooth locus. Since $\Z$ is semi-stable, then so is $\wt{\X}$. 
\qed

\begin{corollary} \label{st-hull.pot}
There exists a finite separable extension $K'/K$
such that $\X_{S'}$, where $S'$ is the normalization of $S$ in
$K'$, is dominated by a semi-stable model of $X_{K'}$.
\end{corollary}

{\noindent\it Proof:} We can suppose that $X$ has semi-stable 
reduction over $S$. Since $X$ has good reduction over an open
dense subset of $S$, we can suppose that $S$ is local. 
By the finiteness Theorem of Grauert-Remmert 
(\cite{GR}, see also \cite{BLR-fr}, Theorem 1.3) applied to
the formal completion of $\X$ along its special fiber, 
there exists a finite Galois extension 
$L/\hat{K}$ such that the normalization of $\X_{\cO_{L}}$
has geometrically reduced special fiber. See \cite{Epp}, p. 247
for how to descend the result to $K$ (Note that 
\cite{Kuh} fills a gap in the proof of a main theorem in \cite{Epp}). 
We then apply 
Proposition \ref{st-hull.ex}. 
\qed 

\begin{remark} \label{DJ}
The corollary  is useful in a recent work of 
Ch. Deninger and A. Werner on vector bundles and representations 
of the fundamental group of $p$-adic curves \cite{DW}. In fact, 
de Jong (\cite{dJ}, Theorem 2.4) already proved it in the 
situation when $S$ is a Noetherian integral excellent scheme 
of any dimension, and when $X$ is an integral curve over $K(S)$. 
The scheme $S'$ is then proper and generically finite over $S$. 
The proof here for one-dimensional $S$ is simpler, and more 
effective in some sense. 
\end{remark} 

\begin{remark} When $S$ is local and complete, the corollary can 
be reformulated in terms of rigid analytic geometry as follows: 
let $\mathcal U$ be a formal covering of $X$. Then after finite 
separable extension of $K$, $\mathcal U$ can be refined to a 
distinguished formal covering $\mathcal V$ with semi-stable 
reduction. As such, the statement can be easily worked out using
\cite{BL}, Theorem 5.5,  and step 2 in the proof of Lemma 7.3, page 377.
The non complete case can then be obtained using 
Proposition \ref{completion}. 
\end{remark} 

\begin{remark} ({\em Effective reduced fiber theorem}). In the case of 
\ref{st-hull.pot}, we can give an effective method to eliminate
the multiplicities of $\X_s$, without using Grauert-Remmert's theorem. 
Suppose that $X$ has semi-stable reduction.
Let $\Gamma$ be an irreducible component of $\X_s$ of geometric
multiplicity $d>1$. If we choose two closed points $P_1, P_2$ of 
$X$ which specialize to two distinct points in the interior of 
$\Gamma$, and if we take $L=K(P_1, P_2)$, then the irreducible components
of $(\X_{\cO_L})'s$ (where $(\X_{\cO_L})'$ denotes the normalization
of $\X_{\cO_L}$) lying above $\Gamma$ are geometrically reduced. 
If $K$ is strictly Henselian, we can bound $[L: K]$ by $d^2$. Note that
$L$ can be chosen to be separable over $K$. If $X$ has not necessarily
semi-stable reduction, then we can bound $[L : K]$ by the max of 
$d^2$ and a constant depending only on $g$. 
\end{remark} 

\begin{emp}{\bf The stable hull}. 
Now let us construct a minimal semi-stable model dominating $\X$. 
Let $\Z$ be a locally complete intersection (e.g., regular or
semi-stable) model of $X$ over $S$. 
Let $\omega_{\Z/S}$ be the (invertible) dualizing sheaf of $\Z/S$.
Recall that a {\em $(-2)$-curve} on $\Z$ is an
irreducible component $\Gamma$ of a closed
fiber $\Z_s$ such that $\deg\omega_{\Z/S}|_{\Gamma}=0$. If $\Z$ is
semi-stable and $k(s)$ is algebraically closed, and $\Gamma$ is not
a connected component of $\Z_s$, then this is equivalent
to $\Gamma\simeq\mathbb P^1_{k(s)}$ and $\Gamma$ meets the other
irreducible components at exactly two points.
Recall that the {\em exceptional locus} of
a birational projective morphism $\pi: \Z\to \X$ is by definition 
the complementary of $\pi^{-1}(U)$, where $U$ is the biggest open 
subscheme of $\X$ such that $\pi^{-1}(U)\to U$ is an isomorphism. 
When $\X$ is normal, the exceptional locus is equal to the union 
of the prime divisors of ${\Z}$ which map to 
closed points in $\X$. A semi-stable model $\Z$ dominating
$\X$ will be called {\em relatively minimal} if 
there is no semi-stable model 
between $\X$ and $\Z$, except $\Z$ itself. 
\end{emp}

\begin{lemma} \label{contract-2}
Let $\Z$ be a semi-stable model of $X$ over $S$. 
Let $V$ be an effective vertical divisor on $\Z$ such that
for all $s\in S$, no connected component of $\Z_s$ is contained in 
$V$. 
\begin{enumerate}[{\hskip 8pt\rm (a)}] 
\item If $\deg\omega_{\Z/S}|_{\Gamma}\le 0$ for all components $\Gamma$
of $V$, then there exists a contraction map $\Z\to \W$ of $V$ 
and $\W$ is semi-stable. 
\item If there exists a contraction map $\Z\to \W$ of $V$ with
$\W$ semi-stable, then $\deg\omega_{\Z/S}|_{\Gamma}\le 0$ for at
least one irreducible component $\Gamma$ of $V$.
\end{enumerate}
\end{lemma} 

{\noindent\it Proof:} (a) is well-known but we were not able to find a proper
reference. We can suppose that $S$ is local. Let $\rho : \wt{\Z}\to \Z$ be the 
minimal desingularization. Then the components $\Theta$ of the
exceptional locus $E$ of $\rho$ are $(-2)$-curves. Let $\wt{V}$ be the
strict transform of $V$ in $\wt{\Z}$. Let us show that $V':=E+\wt{V}$
can be contracted. We have 
$\omega_{\wt{\Z}/S}=\rho^*\omega_{\Z/S}$ because $\Z$ is semi-stable.
Hence $\deg\omega_{\wt{\Z}/S}|_{\Gamma'}\le 0$ for all components $\Gamma'$
of $V'$. If there exists a $\Gamma'$ such that 
$\deg\omega_{\wt{\Z}/S}|_{\Gamma'}< 0$, then $\Gamma'$ is an exceptional
divisor. Let $\pi : \wt{\Z}\to\Z'$ be the contraction of $\Gamma'$. Then
$\omega_{\wt{\Z}/S}=\pi^*\omega_{\Z'/S}(\Gamma')$. We deduce easily that
$\deg\omega_{\Z'/S}|_{\Gamma''}\le 0$ for all $\Gamma''$ in $\pi(V')$. 
So by successively blowing-down exceptional divisors, we can suppose
that $V'$ consists only of $(-2)$-curves. By
Artin's criterion of contractibility (\cite{Liub}, 9.4.7),
$V'$ can be contracted. Therefore $V$ can be contracted. 

It remains to see that $\W$ is semi-stable. 
Let $\cO_{K'}$ be a discrete valuation ring
containing $\cO_S$ and let $S'=\Spec\cO_{K'}$. It is enough to
show that $\W_{S'}$ is semi-stable.  
The map $\Z_{S'}\to \W_{S'}$ is the contraction of $V_{k(s')}$.
We have $\omega_{\Z_{S'}/S'}=\omega_{\Z/S}\otimes\cO_{S'}$. 
If $\Gamma'$ is a component of $\Z_{s'}$ lying over $\Gamma\subseteq V$,
then 
$$[k(s'): k(s)]\deg_{k(s')}\omega_{\Z_{S'}/S'}|_{\Gamma'}=
[k(\Gamma') : k(\Gamma)]\deg_{k(s)} \omega_{\Z/S}|_{\Gamma}\le 0.$$ 
So we can reduce the lemma to the case $k(s)$ algebraically closed. 
Let $\Gamma\subseteq V$. Then 
$p_a(\Gamma)=0$, 
$\Gamma\simeq \mathbb P^1_{k(s)}$, and $\Gamma$ meets the other
components of $\Z_s$ in at most two points. 
Now it is well-known that $\W_s$ is semi-stable 
(see for instance \cite{Liub}, 10.3.31).

(b) The previous computations show that 
$\deg \omega_{\Z/S}|_{\Gamma}=\deg\omega_{\wt{\Z}/S}|_{\wt{\Gamma}}$
for any irreducible component $\Gamma$ of $\Z_s$. Let $\wt{\W}$ be the
minimal desingularization of $\W$. 
Suppose that $\wt{\Z}\to \wt{\W}$ is not an isomorphism. 
Let $\Gamma$ be a component of $\Z_s$ whose strict transform 
in $\wt{\Z}$ is an exceptional divisor contracted into a closed 
point in $\wt{\W}$. Then $\Gamma\subseteq V$ and 
$\deg\omega_{\Z/S}|_{\Gamma}=\deg\omega_{\wt{\Z}/S}|_{\wt{\Gamma}}<0$.
If $\wt{\Z}=\wt{\W}$, then 
$\deg\omega_{\Z/S}|_{\Gamma}=\deg\omega_{\wt{\W}/S}|_{\wt{\Gamma}}=0$ for all
$\Gamma$ in $V$ because $\wt{\Gamma}$ (the strict transform of $\Gamma$)
is a $(-2)$-curve in $\wt{\W}$. 
\qed 

\begin{proposition} \label{st-hull.st}
Let $\X$ be a model of $X$ over $S$ dominated by a semi-stable model. 
\begin{enumerate}[{\hskip 8pt\rm (a)}]
\item A semi-stable model $\Z$ dominating $\X$ is relatively 
minimal if and only if for all irreducible components 
of the exceptional locus of $\Z\to\X$,
we have $\deg\omega_{\Z/S}|_{\Gamma}>0$.
\item Let $\wt{\X}\to\X$ be the minimal desingularization
of $\X$. Let $\wt{\X}\to \W$ be the contraction of the $(-2)$-curves
contained in the exceptional locus of $\wt{\X} \to \X$. 
Then $\W$ is the stable hull of $\X$.
\end{enumerate}
\end{proposition} 

{\noindent\it Proof:} (a) is an immediate consequence of
Lemma \ref{contract-2}. 

(b) By \ref{contract-2}(a), $\W$ is semi-stable and dominates $\X$. 
Let $\Z$ be a semi-stable model dominating $\X$, and 
relatively minimal. Let us first show that $\wt{\X}$ 
dominates $\Z$. Let $\wt{\Z}$ be the 
minimal desingularization of $\Z$. Then $\wt{\Z}$ dominates 
$\wt{\X}$. Suppose that $\wt{\Z}\to\wt{\X}$ is not an isomorphism.
Let $\Theta$ be an exceptional divisor of $\wt{\Z}$ mapped
to a closed point of $\wt{\X}$. Since $\wt{\Z}\to \Z$ is minimal,
$\Theta$ maps to an irreducible component $\Gamma$ of $\Z_s$ which
is then contained in the exceptional locus of $\Z\to\X$. We have
$\deg\omega_{\Z/S}|_{\Gamma}=\deg\omega_{\wt{\Z}/S}|_{\Theta}<0$.
Contradiction. Therefore $\wt{\Z}\simeq \wt{\X}$ and $\wt{\X}\to \Z$
consists in contracting some $(-2)$-curves in $\X_s$. 
Hence $\Z$ dominates $\W$. 
\qed 

\begin{remark} A {\em Du Val model} of $X$ over $S$ is a model $\X$ 
such that $\deg\omega_{\wt{\X}/S}|_{\Gamma}=0$ where $\wt{\X}\to\X$ is
the minimal desingularization and where $\Gamma$ is any irreducible
component of the exceptional locus of $\wt{\X}\to \X$. 
The above results (\ref{contract-2} and \ref{st-hull.st}) 
still hold when ``semi-stable'' is replaced by ``Du Val'', except
the base change property. The point is that Du Val models do not 
commute with base change. 
\end{remark} 

\begin{remark} \label{st-h-aut}
Let $\W$ be semi-stable and dominating $\X$.
Then $\W$ is the stable hull of $\X$ if and only if 
${\rm Aut}_{\X_{\bar{s}}}(\W_{\bar{s}})$ is finite for all $s\in S$.
\end{remark} 

\begin{emp} \label{pf-ssh}
{\bf Proof of Theorem \ref{ss-h}}. (a) is contained in 
Corollary \ref{st-hull.pot} and Proposition \ref{st-hull.st}. 
(b) It is enough to show that $\W_{S'}$ is relatively minimal. 
Let $\Gamma'$ be an irreducible component of $\W_{s'}$ contained in the
exceptional locus of $\W_{S'}\to\X_{S'}$. The image of $\Gamma'$ in
$\W$ is in the exceptional locus of $\W\to\X$. Similarly to the
proof of Lemma \ref{contract-2}(a), we have 
$ \deg \omega_{\W_{S'}/S'}|_{\Gamma'} > 0$. Therefore $\W_{S'}$ is
relatively minimal by \ref{st-hull.st}(a). 
\qed 

\end{emp}

\begin{remark} Suppose that $g(X)\ge 1$ and $S$ is affine.  
Let $\X$ be the minimal regular model of $X$ over $S$ and let 
$\W$ be the stable hull of $\X_{S'}$ over some extension $S'/S$. 
Let $\X'$ be the stable or minimal regular model of $X_{K'}$ over 
$S'$. Then $H^0(\W, \omega_{\W/S'})=H^0(\X', \omega_{\X'/S'})$.
One should be able to recover some arithmetic informations on $\X$
from the sheaf 
$\omega_{\W/S'}\otimes((\omega_{\X/S})^{\vee}\otimes\cO_{S'})$. 
Let us consider the ideal 
$$\cO_{S'}(-\delta):=H^0(\W, \omega_{\W/S'})\otimes 
(H^0(\X, \omega_{\X/S})^{\vee}\otimes\cO_{S'}).$$
For example, if $X$ is an elliptic curve over $K$, then for 
every closed point $s\in S$, we can show that 
$$12{\rm ord}_s(\delta)={\rm ord}_s(\Delta) + a_s{\rm ord}_s(j),$$ 
where $\Delta$ is the minimal discriminant divisor of $X$ over $S$, 
and $a_s=0$ if $X$ has potentially good reduction at $s$, 
$a_s=1$ otherwise.
\end{remark} 

\begin{emp} \label{marked} {\bf Marked curves}. Recall that a 
(proper) {\em marked curve} $Z\to T$ over a scheme $T$ is a 
proper flat scheme of relative dimension $1$ over $T$ 
endowed with a finite set $M\subset Z(T)$ of sections
with pairwise disjoint supports 
contained in the smooth locus of $Z/T$ (for our purpose, it is not 
necessary to order these sections).  
Note that if $T$ is irreducible with generic point $\xi$, 
then $M$ is determined by its generic fiber $M\cap Z_{\xi}$.  
We say that $(Z, M)$ is {\em semi-stable} if $Z \to T$ is semi-stable.
We say that $(Z, M)$ is {\em stable} if it is semi-stable and if
for any geometric point $\bar{t}$ of $T$, $Z_{\bar{t}}$ is connected
and for any irreducible component $\Gamma$ of $Z_{\bar{t}}$, $\Gamma$ 
meets the other components in at least 
$1-(2p_a(\Gamma)-2)-|M\cap \Gamma|$ points. This amounts to 
say that $\omega_{Z/T}(M)$ is ample. 

A {\em morphism} of marked curves $(Z, M) \to (Z', M')$ over $T$ 
is a $T$-morphism $f : Z \to Z'$ such that $f(M)\subseteq M'$. 

Let $(X, M)$ be a smooth marked curve over $K=K(S)$. A {\em marked model}
of $(X, M)$ over $S$ is a marked curve $(\X, \mathcal M)$ over
$S$ whose generic fiber is isomorphic to $(X, M)$. Since $\mathcal M$
is uniquely determined by $M$ and $\X$, we will omit $\mathcal M$
in the notation $(\X, \mathcal M)$ and we will simply say
$\X$ is a marked model of $(X, M)$.  

Let $\X$ be a model of $X$ over $S$. The 
{\em stable marked hull} of $\X$ is the minimal semi-stable 
marked model of $(X, M)$ dominating $\X$. Note that a stable
marked hull is not necessarily a stable marked curve. 
\end{emp} 

\begin{corollary} \label{st-hull.m} Let $(X, M)$ be a smooth
marked curve over $K$ and let $\X$ be a (non-marked) model of $X$. Then 
after finite separable 
extension of $K$, $\X$ admits a stable marked hull. 
More precisely, if $\X$ has a stable hull over some extension $S'/S$, 
then it has a stable marked hull over $S'$. 
Moreover,
the formation of stable marked hull commutes with flat base change
of  Dedekind schemes. 
\end{corollary}

{\noindent\it Proof:} We can suppose that $\X$ has a stable hull $\Z$ over
$S$. Let $\wt{\Z}$ be a desingularization of $\Z$, let 
$\overline{M}$ be the Zariski closure of $M$ in $\wt{\Z}$ and
let $\wt{\Z}_{s_1},...,\wt{\Z}_{s_n}$ be the fibers such that
$\overline{M}\to \overline{M}\cap \wt{\Z}_{s_i}$ is not injective. Let 
$\Z'\to\wt{\Z}$ be an embedded resolution $\overline{M}+\sum_i\wt{\Z}_{s_i}$ 
in $\wt{Z}$ so that the Zariski closure $\mathcal M'$ of $M$ in $\Z'$ is a 
disjoint union of sections (contained in the smooth locus because $\Z'$ is 
regular). Then $\Z'$ is a semi-stable marked model dominating 
$\X$.  

Let $\W$ be any semi-stable marked model of $(X, M)$ over $S$. 
Similarly to the non marked case, we can show that 
$\W$ is relatively minimal if and only if 
$\omega_{\W/S}(\mathcal M)|_{\Gamma}$, where $\mathcal M$ is the
Zariski closure of $M$ in $\W$,  
has positive degree for all $\Gamma$ in the exceptional locus of
$\W\to\X$. This then implies that the stable marked hull is obtained
by contracting prime divisors $\Gamma$ in the exceptional
locus of $\Z'\to\Z$ such that 
$\deg(\omega_{\Z'/S}(\mathcal M')|_{\Gamma}) \le 0$, and that 
the stable marked hull commutes with flat base change. 
\qed 

\begin{remark} \label{stm}
If $(X,M)$ is stable (meaning that $2g(X)-2+|M| \ge 1$) and if 
$X$ has semi-stable reduction over $S$, then
there exists a semi-stable marked model $\X$ of $(X, M)$ over $S$ 
and minimal for this property. This model is the stable marked
model of $(X, M)$. It is characterized by the property that
for all irreducible components $\Gamma$ of $\X_s$, one has 
$\deg\omega_{\X/S}(\mathcal M)|_{\Gamma} > 0$, where $\mathcal M$
is the Zariski closure of $M$ in $\X$. As above, this implies that
the stable marked model commutes with base change. 
\end{remark} 

\end{section} 

\begin{section}{Semi-stable models of finite covers} 

We (re)prove that any finite morphism of projective smooth
curves over $K$ extends, after finite separable extension of $K$,
to a finite morphism of semi-stable models.

\begin{definition} \label{model-f}
Let $f : X \to Y$ be a finite 
morphism of smooth connected projective curves over $K(S)$. 
A {\em model (or extension) of $f$ over $S$} consists in a 
morphism $\X \to \Y$ 
over $S$ extending $f$, where $\X$ and $\Y$ are models over $S$ 
of $X$ and $Y$ respectively. A model of $f$ is said to be 
{\em finite} if it is a finite morphism, and {\em semi-stable} 
(see \cite{Col}) if it is finite and if $\X$, $\Y$ are semi-stable. 
We say that a model $\X \to \Y$ of $f$ 
dominates another one $\X'\to \Y'$ if there are birational morphisms
$\X\to \X'$, $\Y\to\Y'$ making the following diagram commutative
$$\begin{CD}
\X @>>> \Y \\
@VVV @VVV \\
\X' @>>> \Y'
\end{CD}
$$
A model of $f$ is {\em stable} if it is semi-stable and minimal
(for the domination relation) among semi-stable models of $f$. 
If $\X, \Y$ are respective models of $X, Y$. Then the 
semi-stable model $\X'\to\Y'$ of $f$ such that $\X'$ dominates
$\X$, $\Y'$ dominates $\Y$, and which is minimal for these
property, is called the {\em stable hull} of (the rational map)
$\X\dasharrow\Y$.
We can obviously make similar definitions for marked curves. 
\end{definition}

\begin{remark} For a given $f : X\to Y$,
the semi-stable models are not unique: let $\X\to\Y$ be a semi-stable
model of $f$, let $\Y'\to \Y$ be a blowing-up along a closed point,
then the stable hull of $\X\dasharrow\Y'$ (see \ref{stable-e}) 
is a new semi-stable model. 
\end{remark} 

\begin{emp}{\bf Decomposition of inseparable morphisms}. 
Let us first deal with purely inseparable morphisms $X\to Y$. The
next two statements are well-known at least over perfect
base fields. 

\begin{lemma} Let $K$ be a field of characteristic $p>0$. Let
$E/F$ be a finite extension of function fields of one variable
over $K$, with $E$ separable over $K$. Then there exists a unique 
purely inseparable sub-extension $L/F$ of $E/F$ such that $E/L$ 
is separable. Moreover, $F=K L^{p^r}$ for some $r\ge 0$.
\end{lemma} 

\noindent{\it Proof:} Let $F_s$ be the separable closure of $F$ in
$E$. By \cite{Liub}, Cor. 3.2.27 (here we use the hypothesis $E$ separable 
over $K$), there exists $r\ge 0$  such that $F_s = K E^{p^r}$. Let
$$L:=\{ e \in E \ | \ e^{p^r}\in F \}.$$
Then $L/F$ is a purely inseparable extension, $ F= K L^{p^r}$, and
$E/L$ is separable because otherwise $L\subseteq K E^p$ and
$F_s\subseteq K E^{p^{r+1}}$. 
The uniqueness of $L$ is obvious because it is necessarily equal
to the radicial closure of $F$ in $E$. 
\qed

\begin{proposition} \label{decomp}
Let $f : X \to Y$ be a finite morphism of 
normal connected curves over a field $K$ of characteristic $p>0$. 
Suppose that $X$ is smooth. Then $f$ can be decomposed into a 
finite separable morphism $X \to Z$ followed by $Z\to Y$ which
can be identified to a Frobenius map $Z\to Z^{(p^r)}$ for some
$r\ge 0$. Moreover $Z$ is smooth.
\end{proposition} 

{\noindent\it Proof:} Let $L$ be the radicial closure of $K(Y)$ in
$K(X)$, and let $Z$ be the normalization of $Y$ in $L$. Then
$f$ induces a finite separable morphism $X\to Z$. It is 
flat because $Z$ is regular of dimension $1$. Let $\bar{K}$ be
an algebraic closure of ${K}$, then 
$X_{\bar{K}}\to Z_{\bar{K}}$ is flat, hence $Z_{\bar{K}}$ is
regular. Finally, $Z\to Y$ can be identified to 
$Z\to Z^{(p^r)}$ by \cite{Liub}, 7.4.21.
\qed 
\smallskip 

Note that $f$ can also be decomposed into $X\to X^{(p^r)}$
followed by a separable morphism $X^{(p^r)}\to Y$.
\end{emp}

\smallskip 

\begin{emp}{\bf Semi-stable models}. 
\end{emp} 

\begin{lemma} \label{mods} 
Let $f_i : X\to Y_i$, $i=1,...,n$ be finite surjective morphisms
of integral projective varieties over  
$K$ and let $\Y_i$ be a model of
$Y_i$ over $S$. Then there exists a model $\X$ of $X$ over $S$
such that $\X$ dominates $\Y_i$ (that is, $X\to Y_i$ extends to 
$\X\to \Y_i$) for all $i$.
\end{lemma} 

{\noindent\it Proof:} For $N$ big enough, we have a closed
immersion $g : X\to \mathbb P^N_K\times_K (\prod_{K,i} Y_i)$ induced by the
projective morphism $\prod_i f_i : X\to \prod_{K,i} Y_i$. Now take
$\X$ to be the Zariski closure of $g(X)$ in 
$\mathbb P^N_S\times_S (\prod_{S,i} \Y_i)$, endowed with the reduced
structure. Note that if $X$ is geometrically reduced, we can also 
use Lemma \ref{join} with $\X_i=N(\Y_i, K(X))$. 
\qed 

\begin{proposition} \label{semi-stable} 
Let $S$ be a Dedekind scheme, and let $f : X\to Y$ be a finite
morphism of smooth geometrically connected projective curves over $K:=K(S)$.
Let $\X, \Y$ be respective models of $X$ and $Y$. 
Then there exists a finite separable extension $K'/K$ such that
over the normalization $S'$ of $S$ in $K'$, 
the cover $X_{K'}\to Y_{K'}$ extends to a finite morphism
$\X'\to \Y'$, where $\X'$ (resp. $\Y'$) is a semi-stable model of 
$X$ (resp. $Y$) over $S$ dominating $\X_{S'}$ (resp. $\Y_{S'}$). 
\end{proposition} 

{\noindent\it Proof:} Let $X\to Z\to Y$ be the decomposition
given by Proposition \ref{decomp}. Let $\hat{X}\to Z$ be the Galois
closure of $X\to Z$. After a finite separable extension of $K$, 
$\hat{X}$ is smooth over $K$. Let 
$\hat{\X}_0/S$ be a model of $\hat{X}$ dominating $\X$ and 
$\Y$ (\ref{mods}). Let  $G:={\rm Gal}(K(\hat{X})/K(Z))$.
By \ref{equiv}, after a finite separable extension of $K$, 
there exists a semi-stable model $\hat{\X}$ of
$\hat{X}$ endowed with an action of $G$ and dominating 
$\hat{\X}_0$. 
Let $\X'=\hat{\X}/H$ where $H={\rm Gal}(K(\hat{X})/K(X))$, and
$\Z'=\hat{\X}/G$. Then $\X'\to \Z'$ is a finite morphism of semi-stable
models of $X$ and $Z$ respectively (\cite{Ray1}, Prop. 5). 
Let $\Y'=\Z'^{(p^r)}$. Then the canonical map $\Z'\to \Y'$ is finite
and $\Y'$ is semi-stable (loc. cit., or \cite{Liub}, 
Exerc. 10.3.19(a)). Since $\hat{\X}$ dominates $\Y$ and 
is finite over $\Y'$, we see easily that $\Y'$ dominates $\Y$
(use for instance \cite{LL}, 4.1). Hence the proposition is proved
with $f'$ equal to the composition $\X'\to \Z'\to \Y'$. 
\qed 

\begin{remark} If $S$ is any Noetherian integral excellent 
scheme, then using the result of de Jong \cite{dJ} as quoted
in \ref{DJ}, we see that the proposition is still true.
\end{remark}

The next corollary was known for separable morphisms
(\cite{Col} when $K$ is complete; \cite{LL}, Remark 4.6 when $g(X)\ge 1$). 

\begin{corollary} \label{sst.red}
Let $f : X\to Y$ be a finite morphism of smooth geometrically
connected projective curves over $K$. Then after a finite 
separable extension of $K$, there exists a finite morphism 
$\varphi : \X \to \Y$ of semi-stable models of $X, Y$ respectively.
\end{corollary} 

{\noindent\it Proof:} Apply Proposition \ref{semi-stable} to any
pair of models of $X, Y$.
\qed

\begin{lemma} \label{split} 
Let $S$ be local. Let $F$ be a finite closed subset of $X$. Let $\Z$
be a semi-stable model of $X$ over $S$. Then there exists an integer
$d>0$ such that for any finite extension $\cO_{K'}/\cO_S$ of discrete
valuation rings with ramification index divisible by $d$, if 
$\wt{\Z}'$ denotes a desingularization of $\Z_{\cO_{K'}}$,
then the Zariski closure of $F_{K'}$ in $\wt{\Z}'$ is contained
in the smooth locus of $\wt{\Z}'$. 
\end{lemma} 

{\noindent\it Proof:} Let $\alpha\in F$, and let $x\in \Z_s$ be 
a singular specialization of $\alpha$. The local ring of 
an \'etale neighborhood of $x\in \Z$ is 
isomorphic to $\cO_{K}[[u,v]]/(uv-a)$, with $a$ a power of a
uniformizing element of $\cO_K$. Let $u(\alpha)$ be the image of 
$u$ in $K(\alpha)$. Then $|a|<|u(\alpha)|<1$. After an extension
of  big enough ramification index, $|u(\alpha)|$ belongs to
$|K|$. If this condition is satisfied for all $\alpha\in F$ and
for all singular specializations of $\alpha$, then it is easy
to see that the specializations of $F$ in $\wt{\Z}'$ are
smooth points. Indeed, if $\Z$ is regular, the parameter 
$a$ in the above local ring is a uniformizing element, 
hence $|a| < | u(\alpha)| < 1$ cannot hold in $|K|$, so $F$
must specialize to smooth points. 
\qed

\begin{proposition} \label{ss-m}
Let $f : (X, M)\to (Y, N)$ be a finite morphism
of smooth geometrically connected marked projective curves over $K$, let 
$\X, \Y$ be respective models of $X, Y$. Then after a finite
separable extension of $K$, there exists a semi-stable marked
model $\X'$ (resp. $\Y'$) of 
$(X, M)$ (resp. $(Y,N)$) such that $\X'$ and $\Y'$ dominate
respectively $\X$ and $\Y$, and $f$ extends to a finite morphism 
$\X'\to \Y'$. 
\end{proposition} 

{\noindent\it Proof:} After enlarging $K$ and replacing $\X$ and 
$\Y$ be their respective stable marked hull, we can suppose that
$\X$, $\Y$ are semi-stable and that the Zariski closure 
$\overline{M}$ of $M$ in $\X$ is a disjoint union of sections, 
and the same for $N$ in $\Y$. 
Let $X\to Z\to Y$ be the decomposition as
given by Proposition \ref{decomp}, and let $\hat{X}\to Z$ be the
Galois closure of $X\to Z$. 
Let $f : \hat{X}\to X$ and $g : \hat{X}\to Y$ be
the canonical morphisms. 
By lemma \ref{split}, after a finite separable 
extension, and after replacing $\hat{\X}$ by its minimal desingularization
(the group $G$ still acts), we can suppose that the Zariski closure of 
$f^{-1}(M)\cup g^{-1}(N)$ in $\hat{\X}$ is contained in the
smooth locus. Then the Zariski closure of 
$M$ in $\X'$ is contained in the smooth locus because
$\hat{\X}_{\rm sm}/H$ is smooth, and it is a disjoint union
of sections because $\X'$ dominates $\X$ and 
$\overline{M}$ is already a disjoint union of sections. 
Hence $\X'$ is semi-stable marked for $(X, M)$. 
The same arguments hold for $\Y'$. 
\qed 

\end{section} 

\begin{section}{Stable hull of a morphism} 

\begin{definition} Let $f : X \to Y$ be a finite morphism of 
connected smooth projective curves over $K$. Let $\X, \Y$ be 
respective models of $X, Y$ over $S$. The {\em finite hull of 
$\X\dasharrow \Y$} is a finite model $\X^{\rm f}\to \Y^{\rm f}$ of
$f$ over $S$, such that $\X^{\rm f}$ and $\Y^{\rm f}$ are normal 
models of $X, Y$ dominating respectively $\X$ and $\Y$, 
and which is minimal (for the domination relation) 
with respect to these properties. 
\end{definition}

\begin{lemma} \label{normalization} 
Let $\Y$ be an integral scheme locally of finite type over
$S$, let $L$ be a finite extension of $K(\Y)$, separable over
$K=K(S)$, and let $\X$ be the normalization of $\Y$ in $L$. 
Then $\X$ is finite over $\Y$. 
\end{lemma} 

{\noindent\it Proof:} 
The assertion is of course trivial if $S$ is excellent. 
Let $y\in\Y_s$, let $A=\cO_{\Y,y}$, and let $B$ be the integral closure of 
$A$ in $L$. We have to show that $B$ is finite over $A$.
Let $\cO_K=\cO_{S,s}$. Let $C$ be the integral closure of 
${A}\otimes_{\cO_K} \hat{\cO}_K$ in $L\otimes_{K} \hat{K}$.
The latter is reduced (because $L$ is separable over $K$)
and finite over  $K(\Y)\otimes_{K}\hat{K}$, the total ring of fractions
of $A\otimes\hat{\cO}_K$.
Since $\hat{\cO}_K$ is excellent, $C$ (and thus $B\otimes\hat{\cO}_K$)
is finitely generated over $A\otimes\hat{\cO}_K$. Then 
Nakayama's lemma implies that $B\otimes\hat{\cO}_K$ is generated
over $A\otimes\hat{\cO}_K$ by finitely many elements of $B$. These
elements generate $B$ over $A$ because $\cO_K\to\hat{\cO}_K$
is faithfully flat.

Note that the proof still work if $S$ is any Noetherian integral
scheme such that $\hat{\cO}_{S,s}$ is reduced for all $s\in S$. 
\qed

\begin{lemma} \label{join} Let $X$ be an integral projective 
variety over $K$, and let $\X_1, \X_2, ..., \X_n$ 
be models of $X$ over $S$. 
\begin{enumerate} [{\hskip 8pt\rm (a)}]
\item There exists a smallest model $\X$ of $X$ dominating 
$\X_i$ for all $i$. Let us denote $\X$ by $\X_1\vee\cdots\vee\X_n$. 
\item If $X$ is geometrically integral, then for any flat morphism
of Dedekind schemes $S'\to S$, we have 
$(\X_1\vee\cdots\vee\X_n)_{S'}=(\X_1)_{S'}\vee\cdots\vee (\X_n)_{S'}$. 
\item If $\dim X=1$, then
every irreducible component of $(\X_1\vee\cdots\vee\X_n)_s$ 
dominates an irreducible component of $(\X_i)_s$ for some $i$. 
\end{enumerate} 
\end{lemma} 

{\noindent\it Proof:} The proof is similar to that of \ref{equiv-lm}.
Let $\mathcal P$ be the fiber product 
$\prod_{S,i}\X_i$ over $S$. Then the diagonal map makes $X$ 
a closed subscheme of $\mathcal P_K$. Let $\X$ be the Zariski 
closure of $X$ in $\mathcal P$ endowed with the reduced 
structure. Then $\X$ is a model of $X$ over $S$ dominating 
the $\X_i$'s.  
Let $\Z$ be a model of $X$ over $S$ dominating the $\X_i$'s. Then
we have a natural morphism $f : \Z\to \mathcal P$ whose image
$f(\Z)$ is irreducible, with generic fiber $\X_K$. Hence 
$f(\Z)=\X$ and $f$ factorizes through $\Z\to\X\to\mathcal P$. 
So $\X$ is minimal. If $X$ is geometrically integral, then 
$\X_{S'}$ is an integral closed subscheme of 
$\mathcal P':=\prod_{S', i}(\X_i)_{S'}$, with generic fiber isomorphic to the
diagonal of $\mathcal P'_{K'}$. By construction, $\X_{S'}$ is 
equal to $(\X_1)_{S'}\vee\cdots\vee(\X_n)_{S'}$. 

Let $s\in S$. Then $\X_s$ is a closed subscheme of 
$\mathcal P_s=\prod_{k(s),i}(\X_i)_s$, pure of dimension $\dim X$. 
Let $\Gamma$ be an irreducible component of $\X_s$.
If $\dim X>0$, then the image of $\Gamma$ in $\X_{i_0}$ has 
positive dimension for some $i_0$. If $\dim X=1$, then 
the image of $\Gamma$ in $(\X_{i_0})_s$ is an irreducible 
component. 
\qed

\begin{proposition} \label{finite-h} 
Let $\X\dasharrow \Y$ be a rational map. 
Suppose that either $S$ is local and Henselian, or
$\X$ and $\Y$ are semi-stable or regular. Then the 
finite hull of $\X\dasharrow\Y$ exists. It commutes with base 
changes in the following sense : let $S'\to S$ be a flat morphism
of Dedekind schemes. Then the finite hull of $\X_{S'}\dasharrow \Y_{S'}$
exists and is equal to $((\X^{\rm f})_{S'})^{\sim}\to 
((\Y^{\rm f})_{S'})^{\sim}$, where $^{\sim}$ means normalization.
\end{proposition} 

{\noindent\it Proof:} The rational map $\X\dasharrow\Y$ is defined and
finite above an open dense subset of $S$. So we can suppose that 
$S=\Spec\cO_K$ is local with closed point $s$. The case when $\X, \Y$
are semi-stable or regular is easily reduced to the Henselian case by 
passing to the completion of $\cO_{K}$ and using Proposition \ref{completion}. 
Suppose that $\cO_K$ is Henselian. Let $\X_1$ be the normalization of
$\X\vee N(\Y, K(X))$ (see \ref{join}). We have a morphism $\X_1\to \Y$. 
By \cite{LL}, Lemma 4.14, 
there exists a (unique) normal model $\Y^{\rm f}$ such that the 
rational map $\X_1\dasharrow \Y^{\rm f}$ is quasi-finite and 
surjective in codimension $1$ 
(in other words : if $\mathcal U$ is the domain of definition of
$\X_1\dasharrow\Y^{\rm f}$, then $\mathcal U_s$ is dense in $(\X_1)_s$, 
$\mathcal U_s\to \Y_s$ is quasi-finite and has dense image). 
Since $\X_1\to\Y$ is a morphism, $\Y^{\rm f}$ dominates $\Y$. 
Let $\X^{\rm f}$ be the normalization of $\Y^{\rm f}$ in $K(X)$. 
Then it is easy to see that $ \X^{\rm f}\to \Y^{\rm f}$ is the finite
hull of $\X\dasharrow \Y$ (use \cite{LL}, 4.1 for instance).

It remains to prove the base change property. 
The rational map $(\X_{S'})^{\sim} 
\dasharrow ((\Y^{\rm f})_{S'})^{\sim}$ 
is quasi-finite and surjective in codimension $1$. 
By the above construction, $(\Y_{S'})^{\rm f}=((\Y^{\rm f})_{S'})^{\sim}$
and $(\X_{S'})^{\rm f}$ is the normalization of 
$((\Y^{\rm f})_{S'})^{\sim}$ in $K(X_{K(S')})$. 
Since $((\X^{\rm f})_{S'})^{\sim}\to ((\Y^{\rm f})_{S'})^{\sim}$
is finite, $(\X_{S'})^{\rm f}$ is isomorphic to 
$((\X^{\rm f})_{S'})^{\sim}$. 
\qed

\begin{thm} \label{stable-e}
Let $S$ be a connected Noetherian regular scheme of dimension $1$, 
let $f : X \to Y$ be a finite morphism of smooth geometrically 
connected projective 
curves over $K:=K(S)$. Let $\X, \Y$ be respective models of 
$X, Y$ over $S$. Then after a finite separable extension of $K$, 
$\X\dasharrow\Y$ admits a stable hull $\X'\to \Y'$. Moreover, 
the formation of $\X'\to \Y'$ commutes with flat base change.
\end{thm}

{\noindent\it Proof:} After a finite separable extension
of $K$, we can suppose that there exists a semi-stable model 
$\X_{\infty}\to \Y_{\infty}$ of $f$ dominating $\X\dasharrow \Y$ 
(Cor. \ref{sst.red}). Let us show that $\X\dasharrow\Y$ then admits
a stable hull over $S$. 
Consider the stable hull 
$\Y_1$ of $\Y$ and the stable hull $\X_1$ of $N(\Y_1, K(X))$. Then
$\X_1$ and $\Y_1$ are dominated by $\X_{\infty}$ and $\Y_{\infty}$ 
respectively. Let $\X_2\to \Y_2$ be the finite hull of $\X_1\to \Y_1$. 
Then it is also dominated by $\X_{\infty}\to \Y_{\infty}$. Now restart 
again the process of taking stable hull and finite hull with 
$\X_2\to \Y_2$. 
We construct in this way an increasing sequence of (normal) models 
$\X_n\to \Y_n$ of $f$ over $S$ 
which are dominated by $\X_{\infty}\to \Y_{\infty}$. This sequence 
is stationary at some rank $n_0$. Then $\X_{n_0}\to \Y_{n_0}$ is 
a semi-stable model of $f$. Note that the construction of 
$\X_n\to \Y_n$ does not depend on the choice of $\X_{\infty}\to \Y_{\infty}$. 
In particular, $\X_n\to \Y_n$
is dominated by any semi-stable model of $f$ dominating $\X\dasharrow \Y$. 
Therefore, $\X_{n_0}\to \Y_{n_0}$ is the stable hull of $\X\dasharrow \Y$. 
Finally, the formation of the stable hull commutes with flat
base change because the stable hull of a model and the finite
hull of a morphism commute with flat base change
(\ref{ss-h}, \ref{finite-h}). 
\qed

\begin{corollary} \label{stable}
Suppose that either $g(X)\ge 2$, or $g(X)=1$ and $X$ has 
potentially good reduction. Then there exists a finite separable 
extension of $K'$ of $K$ such that $X_{K'}\to Y_{K'}$ admits a stable model 
$\X'\to\Y'$ over $S'$, where $S'$ is the integral closure of $S$
in $K'$. Moreover, for any Dedekind scheme $T$ dominating $S'$, 
$\X'_{T}\to \Y'_{T}$ is the stable model of $X_{K(T)}\to Y_{K(T)}$. 
\end{corollary} 

{\noindent\it Proof:} We can suppose that $X$ has semi-stable reduction
over $S$. The cover $X\to Y$ extends to a finite morphism of smooth 
projective models over a dense open subset of $S$. So we can suppose 
that $S=\Spec\cO_K$ is local. Let $\X^{\rm st}$ be the stable 
(resp. smooth projective) model of $X$ if $g(X)\ge 2$ (resp. if $g(X)=1$). 
Suppose first that $\cO_K$ is complete (hence Henselian). Let 
$\X^{\rm st}\dasharrow \Y''$ be the rational map extending 
$X\to Y$ and which is quasi-finite and surjective in codimension $1$ 
(see \cite{LL}, 4.14). Then the stable hull $\X\to \Y$ 
of $\X^{\rm st}\dasharrow \Y''$ is clearly the stable model of 
$X\to Y$. The construction of $\X\to \Y$
commutes with flat base change because that of $\Y''$ and 
the stable hull commute with flat base change. 
If $\cO_K$ is non necessarily complete, we can use Corollary \ref{st-d}.
\qed 

\begin{remark} Let $X\to Y$ be as above. If $X\to Y$ has a semi-stable
model over $S$, then it has a stable model over $S$.
This can be seen in the proof of \ref{stable}. 
If $X\to Y$ is moreover Galois of group $G$,  and if $\X$ is 
the stable model (or smooth model if $g(X)=1$) of $X$ over $S$.
Then the stable model of $X\to Y$ is equal to $\X\to \X/G$.
\end{remark} 

\begin{remark} \label{stable-g}
Suppose moreover that $X\to Y$ is separable, and 
that the Galois closure $\hat{X}$ of $X\to Y$ is smooth and 
geometrically connected over $K$ (which is true after a 
finite separable extension of $K$). 
Let $\hat{\X}$ be the stable model of $\hat{X}$, and let 
$G={\rm Gal}(K(\hat{X})/K(Y))$, $H={\rm Gal}(K(\hat{X})/K(X))$
as in the proof of \ref{semi-stable}. Then we can ask whether
$\hat{\X}/H\to \hat{\X}/G$ is the stable model of $X\to Y$. 
The answer is no in general. Let us give an example with
$X$ and $Y$ having good reduction. 

\medskip
\centerline{\BoxedEPSF{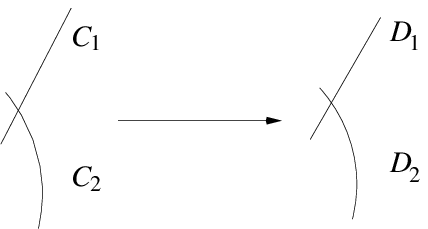}}
\medskip
 
Let $S$ be local, complete, with algebraically closed residue field
$k$. Let $C_1\to D_1$ be a finite separable morphism of degree
$d \ge 3$ with $C_1, D_1\simeq \mathbb P^1_k$, totally ramified 
above some point $y_1\in D_1$, and such that the Galois closure 
$E$ of $C_1\to D_1$ is a curve of genus $g(E)\ge 1$. Let 
$C_2\to D_2$ be a finite separable morphism of degree $d$ 
of smooth connected projective curves over $k$, totally 
ramified above a point $y_2\in D_2$ and such that $g(C_2)\ge 1$. 
Let $D$ be the semi-stable curve over $k$ obtained by identifying
$y_1$ and $y_2$. Let $C$ be the semi-stable curve defined in a
similar way. Then we have a finite morphism $\rho : C\to D$
which is generically \'etale, and such that $C_1\to D_1$, $C_2\to D_2$
have the same ramification index at $y_1$ and $y_2$. 
By \cite{LiuT}, Prop. 5.4, the cover $C\to D$
lifts to a finite morphism $\mathcal C\to \mathcal D$
over $S$ with smooth generic fibers $X, Y$. 
Let $\mathcal C\to\mathcal C_2$ (resp.  $\mathcal D\to\mathcal D_2$)
be the contraction of $C_1$ (resp. of $D_1$). Then $\mathcal C_2$,
$\mathcal D_2$ are smooth, and the canonical morphism 
$\mathcal C_2\to\mathcal D_2$ is the stable model of $X\to Y$. 
Let $\Z$ be the normalization of $\mathcal D$ in $K(\hat{X})$. 
Let $\Theta$ be an irreducible component of $\Z_s$ lying over $D_1$.
Then the separable closure of $k(D_1)$ in $k(\Theta)$ is Galois 
over $k(D_1)$ (\cite{Ser}, I, \S 7, Prop. 20) and contains
$k(C_1)$. Thus $k(\Theta)$ contains a subfield
isomorphic to $k(E)$. In particular, $p_a(\Theta)\ge 1$. 
If $\mathcal C_2\to \mathcal D_2$ is 
equal to $\hat{\X}/H\to \hat{\X}/G$, then $\Z$ dominates $\hat{\X}$
because $\mathcal D$ dominates $\mathcal D_2$.
But $\Theta$ maps to a closed point of $\hat{\X}_s$, thus
$\hat{\X}$ can not be semi-stable. Contradiction. 
\end{remark} 

\begin{remark}
If $X$ has genus $1$ and multiplicative reduction at some point of $S$, 
or if $g(X)=0$, then over any ramified extension of $S$, 
there is no stable model of the identity map $X\to X$.  
The reason we take a ramified extension $S'/S$ is, when $g(X)=1$, 
$X_{K'}$ has no minimal semi-stable model (see the proof of
Lemma \ref{rm-d}).
\end{remark} 

Let us give a characterization of the stable model.

\begin{definition} Let $X_1\to Y_1$, $X_2\to Y_2$ be morphisms 
of schemes over some base scheme $T$. An {\em isomorphism}
of the pairs $(X_1\to Y_1) \to (X_2\to Y_2)$ is given
by $T$-isomorphisms $X_1\to X_2$, $Y_1\to Y_2$ such that the diagram
$$
\begin{CD}
X_1 @>>> X_2 \\
@VVV    @VVV \\
Y_1 @>>> Y_2 \\
\end{CD}
$$
is commutative. We denote by ${\rm Isom}_T((X_1\to Y_1), (X_2\to Y_2))$
the set of these isomorphisms. Now
${\rm Aut}_T(X\to Y)$ has an obvious meaning.
\end{definition}

\begin{proposition} \label{aut-stable}
Keep the hypothesis of {\rm \ref{stable-e}} and
suppose that $g(X)\ge 2$. Let $\X\to\Y$ be a semi-stable model of 
$X\to Y$. Consider the following properties. 
\begin{enumerate} [{\hskip 8pt\rm (i)}]
\item $\X\to\Y$ is the stable model of $X\to Y$ over $S$; 
\item Let $\Gamma$ be any irreducible component of $\Y_s$ such that
$\deg\omega_{\Y_s/k(s)}|_{\Gamma} \le 0$, then there exists an irreducible
component $\Theta$ of $\X_s$ dominating $\Gamma$ 
such that $\deg\omega_{\X_s/k(s)}|_{\Theta}>0$;
\item ${\rm Aut}_{k(\bar{s})}(\X_{\bar{s}}\to\Y_{\bar{s}})$ is finite for
all $s\in S$.
\end{enumerate}
Then {\rm (i)} $\Longleftrightarrow$ {\rm (ii)} $\Longrightarrow$ {\rm (iii)}.
\end{proposition}  

\noindent{\it Proof:} Looking in the proofs of 
\ref{stable-e} and \ref{stable}, we see that $X\to Y$ 
admits a stable model over $S$. Thus $\X\to\Y$ is stable if and 
only if it is a relatively minimal semi-stable model. 
Hence the equivalence  (i) $\Longleftrightarrow$ (ii) is an
immediate consequence of Lemma \ref{contract-2}. 

Suppose that Condition (ii) is satisfied. Then the same condition holds
over $k(\bar{s})$ (with dualizing sheaves on $\X_{\bar{s}}$ and
$\Y_{\bar{s}}$). Consider the natural inclusion 
$${\rm Aut}_{k(\bar{s})}(\X_{\bar{s}}\to\Y_{\bar{s}})\subseteq 
{\rm Aut}_{k(\bar{s})}(\X_{\bar{s}})\times 
{\rm Aut}_{k(\bar{s})}(\Y_{\bar{s}}).$$ 
Note that the right-hand side is not a finite group in general. 
Let $G$ be the subgroup (of finite index) of 
${\rm Aut}_{k(\bar{s})}(\X_{\bar{s}})$
consisting in automorphisms which fix globally each
irreducible component, and let $H$ be the 
subgroup of ${\rm Aut}_{k(\bar{s})}(\Y_{\bar{s}})$
consisting in automorphisms $\tau$ such that $\tau|_{\Gamma}={\rm Id}$
for every irreducible component $\Gamma$ with 
$\deg\omega_{\Y_{k(\bar{s})/k(\bar{s})}}|_{\Gamma}>0$.
Then $H$ is of finite index in ${\rm  Aut}_{k(\bar{s})}(\Y_{\bar{s}})$
because 
$ \{ \tau \in\Aut_{k(\bar{s})}(\Y_{k(\bar{s})}) \  | \
\tau(\Gamma)=\Gamma \}$ is finite for any such $\Gamma$. 
Thus it is enough to show that 
$G':={\rm Aut}_{k(\bar{s})}(\X_{\bar{s}}\to\Y_{\bar{s}})\cap (G\times H)$ is
finite. Let $I$ be the set of irreducible components $\Theta$
of $\X_{\bar{s}}$ such that 
$\deg\omega_{\X_{\bar{s}}/k(\bar{s})}|_{\Theta} >0$. Then
$G|_{\Theta}$ is finite for all $\Theta\in I$. 
Let $(\sigma, \tau)$ be an element of 
${\rm Ker}(G'\to G\times H\to \prod_{\Theta\in I} G|_{\Theta})$.
Condition (ii) implies that $\tau={\rm Id}$ on $\Y_{\bar{s}}$,
hence $\sigma\in {\rm Aut}_{\Y_{\bar{s}}}(\X_{\bar{s}})$. The
latter is a finite group because $\X_{\bar{s}}\to\Y_{\bar{s}}$ is a 
finite morphism. Since the projection map 
${\rm Aut}_{k(\bar{s})}(\X_{\bar{s}}\to\Y_{\bar{s}})\to
{\rm Aut}_{k(\bar{s})}(\X_{\bar{s}}) $ is injective, the above kernel is finite. Hence $G'$ is finite. 
\qed 
\smallskip

\begin{remark}
In general, (iii) does not imply (i). Let $S$ be local, complete
with algebraically closed residue field $k$. 
Let $\pi : \mathbb A^1_k\to \mathbb A^1_k$ be a finite separable cover
with trivial automorphisms group ${\rm Aut}_k(\pi)$. 
Then $\pi$ extends to a finite cover $\mathbb P^1_k\to \mathbb P^1_k$
totally ramified at $\infty$. We can glue $\pi$ with a finite separable
cover $C_2\to D_2$ of smooth projective curves of genus $\ge 2$ over 
$k$ and obtain a finite cover $C\to D$ with finite automorphisms
group (see the construction in \ref{stable-g}) which lifts to a 
finite morphism of semi-stable curves $\mathcal C\to\mathcal D$
over $S$. By the equivalence of (i) and (ii), $\mathcal C\to\mathcal D$
is not stable.
\end{remark} 

\begin{remark} We can see that Theorem \ref{stable-e} holds for 
finite morphisms of 
smooth projective marked curves $(X, M)\to (Y, N)$ : 
{\em let $\X, \Y$ be 
respective models of $X, Y$, then after finite separable
extension, there exists a stable marked hull $\X'\to
\Y'$ of $\X\dasharrow\Y$. Moreover, the formation 
of $\X'\to \Y'$ commutes with flat base change}. The proof is
the same as for \ref{stable-e}, 
except that we replace stable hull of a model by its
stable marked hull, and we use \ref{ss-m} instead of \ref{sst.red}. 
\end{remark}

The next lemma is a generalization of \cite{LL}, 4.4(a) (take 
$N=\emptyset$).
 
\begin{lemma} 
Let $f : (X, M)\to (Y, N)$ be a finite morphism of connected 
smooth projective marked curves over $K$. Suppose that 
$2g(Y)-2 + |N| > 1$ and 
$f^{-1}(N)=M$. If $(X, M)$ has a stable marked model $\X$ over $S$, then
$(Y, N)$ has a stable marked model $\Y$, and $\X\dasharrow \Y$ is
a morphism. 
\end{lemma} 

{\noindent\it Proof:} The stable marked model exists over $S$
when $X$ has semi-stable reduction over $S$. So the existence of
$\X$ implies that of $\Y$ (see \cite{LL}, 4.8 when $g(Y)\ge 1$,
the case $g(Y)=0$ is trivial). It remains to show that the rational
map $\X\dasharrow\Y$ is defined everywhere. Since the stable marked
model commutes with flat base change (\ref{stm}), we can suppose
that $S$ is local with algebraically closed residue field and 
that $\W:=N(\Y, K(X))$ has a stable marked hull $\Z$ over $S$. 
By definition, $\Z$ dominates $\X$. We are going to show that $\Z\to \X$
is an isomorphism. Or equivalently, that $\Z$ is the stable marked
hull of $\X$. Let $\overline{M}$ denote the Zariski closure of $M$ in $\Z$.
It is enough to show that 
$\deg\omega_{\Z/S}(\overline{M})|_{\Theta} >0$ for all irreducible component 
$\Theta$ of $\Z_s$. If $\Theta$ is in the exceptional locus
of $\Z\to\W$, then this is true because $\Z\to\W$ is the stable
marked hull. Suppose that $\Theta$ is the strict transform of 
some irreducible component $\Gamma$ of $\W_s$. 
Let $\Delta$ be the image of $\Gamma$ in $\Y_s$. 
Then every point $y\in \Delta$ which 
is either a singular point of $\Y_s$ or a specialization of $N$
lifts to a point of $\Theta$ which is either a singular point
or a specialization of $f^{-1}(N)=M$ (\cite{LL}, 4.4(b) for singular
points; use going-down property as in \cite{LL}, 4.3 for specializations
of $N$). Since $\Y$ is stable marked, this implies that $\Theta$ contains 
at least three points of $(\Z_s)_{\rm sing}\cup \overline{M}_s$,
hence $\deg\omega_{\Z/S}(\overline{M})|_{\Theta} > 0$. 
\qed

\begin{remark} It is known that in general, $f : X\to Y$ does
not extend to a morphism of the respective minimal regular
models of $X$ and $Y$ over $S$ (see e.g. \cite{CES}, p. 333). However, if
$g(Y)>0$, then $f$ extends to a morphism of the respective
minimal regular models with normal crossings, at least if $k(s)$
is perfect for all closed points $s\in S$. 
\end{remark}

\begin{corollary} \label{stable-m} 
Let $f : (X, M)\to (Y, N)$ be a finite morphism of geometrically
connected smooth projective marked curves over $K$. Suppose that 
$2g(Y)-2 + |N| > 1$ and $f^{-1}(N)=M$. Then, after a finite 
separable extension of $K$, $f$ admits a stable marked model 
$\X\to\Y$ over $S$. This construction commutes with flat base change. 
\end{corollary} 

\begin{remark} \label{mk-ram}
Let $f : X\to Y$ be a finite {\em separable} 
morphism of smooth projective curves. A natural way to mark 
$X$ and $Y$ is to take $M$ equal to the ramification locus of $f$ 
and $N$ equal to the branch locus. By definition $f^{-1}(N)=M$. 
Of course, in general $M$ is not contained in $X(K)$. But if 
$f$ is tamely ramified (e.g., if ${\rm char}(K)=0$), then 
this becomes true over a finite separable extension of $K$
(\cite{LL}, 3.3). Moreover, if $g(Y)\ge 2$, or $g(Y)=1$ and
$f$ is not \'etale, or if $g(Y)=0$ and $g(X)\ge 1$, then 
$2g(Y)-2+|N| > 1$. 
So after again a finite separable extension of $K$, we 
have a canonical way to define a minimal semi-stable reduction of
$X\to Y$ in which the (horizontal) ramification and branch loci 
are finite unions of sections contained in the smooth locus.  
If $(\deg f) !$ is invertible in $\cO_S$, then 
$\X$ and $\Y$ are the respective stable marked models of 
$X$ and $Y$ (\cite{Moc}, \S 3.11, second Lemma).
\end{remark} 

\begin{remark} Let $\hat{X}$, $G$, $H$ be as in Remark \ref{stable-g}. 
Let $\hat{X}$ be marked with its ramification locus over $X$
(rational over $K$ after a finite extension of $K$). 
Let $\hat{\X}$ be its stable marked model. If $G$ has order 
prime to ${\rm char}(k)$, then the stable marked model of $X\to Y$ 
is equal to $\hat{\X}/H\to \hat{\X}/G$ because both sides are
respectively stable marked models of $X$ and $Y$ (see \ref{mk-ram}). 
But this is false if ${\rm char}(k)$ divides $|G|$. Let us
go back to the example of \ref{stable-g} and let $C_1\to D_1$ be moreover 
\'etale outside of $y_1$. 
Let $\mathcal C'\to \mathcal C$, $\mathcal D'\to \mathcal D$ 
be the stable marked hulls of $\mathcal C$, $\mathcal D$ (marked
with horizontal ramification/branch loci). Let $\mathcal C'\to \mathcal C''$,
$\mathcal D'\to \mathcal D''$ be the contraction of 
(the strict transforms of) $C_1$, $D_1$ respectively. Then
it is easy to see that the stable marked model (\ref{stable-m}) of $X\to Y$
is $\mathcal C''\to\mathcal D''$. Similarly to \ref{stable-g}, we
see that it is different from  $\hat{\X}/H\to \hat{\X}/G$. 

It remains to find a cover $C_1\to D_1$ as above. 
Consider the cover\footnote{This example is given by Michel Matignon.} 
defined by the extension $k(D_1) = k(u) \to 
k(C_1)=k(u,v)$, 
with $v^{p^2+p}+v=u$, where $p={\rm char}(k)>0$. Then $C_1\to D_1$ is
\'etale outside of the pole $y_1$ of $u$. Let us show that the
Galois closure $k(E)$ has positive genus. Let $t\in k(E)$ be such 
that $t^{p^2+p}+t=u$ and $t\ne v$. Let $w=t-v$. Then $w$ satisfies
the equation $(w^{p+1}-wt^p-w^pt)^p-w=0$. Hence $w=z^p$ where 
$z=w^{p+1}-wt^p-w^pt$. We have $z=z^{p^2+p}-z^pt^p-z^{p^2}t$, so
$$ (t/z^p)^p+(t/z^p)=1-(1/z)^{p^2+p-1}. $$
This equation defines a $p$-cyclic cover $E'\to \mathbb P^1_k$, with 
conductor $m=p^2+p-1$ at $z=0$, and \'etale elsewhere. Hence
$g(E)\ge g(E')=(p-1)(m-1)/2 \ge 2$. 
\end{remark} 

\begin{remark} If we restrict ourself to the category of regular
models of $X$, then Theorem \ref{ss-h} still holds. More precisely,
given any model $\X$, there exists (after finite separable extension
of $K$) a unique semi-stable and regular model dominating $\X$
and minimal for this property. This model is just the minimal 
desingularization of the stable hull of $\X$. However it does not commute
with base change except when the normalization of $\X$ is smooth. 
On the other hand, Theorem \ref{stable-e}
is no longer true in the setting of regular models. In general, given 
a morphism of models $\X\to\Y$, there is no finite morphisms of
regular models dominating $\X\to\Y$, even after any finite
extension of $K$ (see \cite{LL}, 6.5). 
\end{remark}

\end{section}

\end{document}